\input amstex
\magnification=1200
\documentstyle{amsppt}
\NoBlackBoxes
\def\Z{\Bbb Z}

\def\Cq{\Bbb C_q}

\def\a{\alpha}
\def\l{\lambda}
\def\ph{\varphi}
\def\de{\delta}

\def\Hom{\operatorname{Hom}}
\def\End{\operatorname{End}}
\def\Tr{\operatorname{Tr}}
\def\dimq{\operatorname{dim}_q}
\def\g{{\frak g}}
\def\h{\frak h}
\def\sltwo{{\frak s\frak l _2 }}
\def\sln{\frak{sl}_n}
\def\<{\langle}
\def\>{\rangle}
\def\Uq{{U_q \frak g}}
\def\Uh{{U_q \frak h}}
\def\Un{U_q \frak s\frak l _n }
\def\U2{U_q \frak s\frak l _2 }
\def\o{\otimes}
\topmatter
\title On Cherednik-Macdonald--Mehta identities
\endtitle
\author Pavel Etingof and Alexander Kirillov, Jr
\endauthor
\address
Dept. of Mathematics, Harvard University, Cambridge, MA 02138, USA
\endaddress
\email etingof\@math.harvard.edu\endemail
\address
 Dept. of Mathematics, MIT, Cambridge, MA 02139, USA
\endaddress
\email kirillov\@math.mit.edu  \endemail

\endtopmatter
\document
\head Introduction
\endhead

In this note we give a proof of Cherednik's generalization of
Macdonald--Mehta identities for the root system $A_{n-1}$ using the
representation theory of quantum groups. These identities, suggested
and proved in \cite{Ch2}, give an explicit formula for the integral of
a product of Macdonald polynomials with respect to a ``difference
analogue of the Gaussian measure''. They can be written for any
reduced root system, or, equivalently, for any semisimple complex Lie
algebra $\g$. Assuming for simplicity that $\g$ is simple and
simply-laced, these identities are given by the following formula:
$$
\aligned
\frac{1}{|W|}\int \de_k \overline{\de_k} P_\l\overline{P_\mu}\gamma\, dx
=&q^{\l^2 +(\mu, \mu+2k\rho)} P_\mu(q^{-2(\l+k\rho)})\\
&\quad \times q^{-2k(k-1)|R_+|}\prod_{\a\in R_+}\prod_{i=0}^{k-1}
				(1-q^{2(\a, \l+k\rho)+2i})
\endaligned
\tag 1
$$
where $\l, \mu$ are dominant integral weights, 
$k$ is a positive integer, $P_\l$ are Macdonald
polynomials associated with the corresponding root system, with
parameters $q^2, t=q^{2k}$(see
\cite{M1, M2} or a review in \cite{K2}), $\de_k$ is the $q$-analogue
of $k$-th power of the Weyl denominator $\de=\de_1$:
$$
\de_k=\prod_{\a\in R^+}\prod_{i=0}^{k-1}(e^{\a/2}-q^{-2i}e^{-\a/2}),
\tag 2
$$
and $\gamma$ is the Gaussian, which we define by
$$
\gamma=\sum_{\l\in P} e^\l q^{\l^2},
\tag 3
$$
where $P$ is the weight lattice. We consider
$\gamma$ as a formal series in $q$ with coefficients from the group
algebra of the weight lattice. In a more standard terminology $\gamma$
is called the theta-function of the lattice $P$. All other notations,
which are more or less standard, will be explained below.

These identities were  formulated in the form we use in a
paper of Cherednik \cite{Ch2}, who also proved them using double
affine Hecke algebras (note: our notations are somewhat different from
Cherednik's ones). We refer the reader to \cite{Ch2} for the
discussion of the history of these identities and their role in
difference harmonic analysis. 

As an important intermediate step, we also prove the following
identity for the Gaussian:
$$
\gamma=\biggl(\prod_{\a\in R_+}(1-q^{2(\a, \rho)})
	\biggr)
	\sum_{\nu\in P_+}q^{(\nu, \nu+2\rho)}(\dimq L_\nu) \chi_\nu.
\tag 4
$$
Here $\chi_\nu$ is the character of the  irreducible
finite-dimensional module $L_\nu$ over $\g$, 
and $\dimq L_\nu:=\chi_\nu(q^{2\rho})$ is the quantum dimension 
of $L_\nu$. This identity was
known to experts and  is not difficult to prove; however, we were unable to
locate a proof in the literature.

\subhead{Notations}\endsubhead
 We use the same notations as in \cite{EK1, EK2} with the following
exceptions: we replace $q$ by $q^{-1}$(note that this does not change
the Macdonald's polynomials) and we use the notation $\ph_\l$ for
``generalized characters'' (see below), reserving the notation
$\chi_\l$ for usual (Weyl) characters. In particular, we define
$\overline{e^\l}=e^{-\l}, \bar q=q$, and for $f\in \Cq[P]$, we define
$f(q^\l), \l\in P$ by $e^\mu(q^\l)=q^{(\mu, \l)}$. For brevity, we
also write $\l^2$ for $(\l, \l)$. Finally, we denote by $\int\
dx:\Cq[P]\to\Cq$ the functional of taking the constant term: $\int
e^\l\, dx=\de_{\l, 0}$.

\head 1. Rewriting the Gaussian
\endhead
In this Section, we prove formula (4) for an arbitrary simple Lie algebra
$\g$. 

\proclaim{Proposition 1} Let $\gamma$ be defined by \rom{(3)}. Then
$$
\gamma=\biggl(\prod_{\a\in R_+}(1-q^{2(\a, \rho)})
	\biggr)
	\sum_{\nu\in P_+}q^{(\nu, \nu+2\rho)}(\dimq L_\nu) \chi_\nu.
$$
\endproclaim
\demo{Proof} The proof is straightforward and uses Weyl character
formula along with the following result: if we extend the
definition of $\chi_\nu$ to all $\nu\in P$ by letting 
$\chi_\nu=(\sum (-1)^{|w|}e^{w(\nu+\rho)})/\de$ (recall that
$\de=\de_1$ is the Weyl denominator) then
$\chi_{w.\nu}=(-1)^{|w|}\chi_\nu$, where $w.\nu=w(\nu+\rho)-\rho$. In
particular, $\chi_\nu=0$ if $\nu$ lies on one of the walls, i.e. if
$s_\a.\nu=\nu$ for some root $\a$. The same applies to $\dimq
L_\nu=\chi_\nu(q^{2\rho})$. Using this, we rewrite the
right-hand side of (4) as follows:
$$\aligned
\sum_{\nu\in P_+}q^{(\nu, \nu+2\rho)}(\dimq L_\nu) \chi_\nu=&
    q^{-\rho^2}\sum_{\nu\in P_+}q^{(\nu+\rho)^2}(\dimq L_\nu) \chi_\nu\\
       =&\frac{q^{-\rho^2}}{|W|}\sum_{\nu\in P}q^{\nu^2}
			(\dimq L_{\nu-\rho})  \chi_{\nu-\rho}.
\endaligned
$$
By Weyl character formula, we can write 
$$
(\dimq L_{\nu-\rho}) \chi_{\nu-\rho}
=\frac{1}{\de(q^{2\rho})}\frac{1}{\de} 
	\sum_{w_1, w_2\in W} (-1)^{|w_1w_2|} q^{2(\nu,w_1(\rho))}
	e^{w_2(\nu)}. 
$$
Substituting this in the previous identity, we get 
$$\align
\sum_{\nu\in P_+}q^{(\nu, \nu+2\rho)}&(\dimq L_\nu) \chi_\nu
=\frac{q^{-2\rho^2}}{|W|\de(q^{2\rho})}\frac{1}{\de}
	\sum \Sb w_1, w_2\in W\\ \nu\in P\endSb
		q^{(\nu+w_1(\rho))^2}(-1)^{|w_1w_2|}e^{w_2(\nu)}\\
      =&\frac{q^{-2\rho^2}}{|W|\de(q^{2\rho})}\frac{1}{\de}
	 \sum \Sb w_1, w_2\in W\\ \l\in P\endSb
	   q^{\l^2}(-1)^{|w_1w_2|}e^{w_2(\l-w_1(\rho))}\\
  =&\frac{q^{-2\rho^2}(-1)^{|R_+|}}{|W|\de(q^{2\rho})}
	\sum \Sb w_2\in W\\ \l\in P\endSb
	q^{\l^2}e^{w_2(\l)}
   =\frac{q^{-2\rho^2}(-1)^{|R_+|}}{\de(q^{2\rho})}
	\sum_{\l\in P} q^{\l^2}e^{\l}
\endalign
$$
Simplifying this, we get the statement of the Proposition.
\qed\enddemo

\example{Example} Let $\g=\sltwo$. Then Proposition~1 gives the
following identity, which can be verified directly: 
$$
\sum_{n\ge 0}q^{n(n+2)/2}[n+1](x^n+x^{n-2}+\dots+x^{-n})
=\frac{1}{1-q^2} \sum_{l\in \Z} x^lq^{l^2/2},
\tag 5
$$
where $[n]=\frac{q^n-q^{-n}}{q-q^{-1}}$.
\endexample

 Finally, we note that the Gaussian can be defined in any semisimple
ribbon category $\Cal C$, i.e. a tensor category, with possibly non-trivial
commutativity isomorphism, and a ``Casimir element'' (also called
``ribbon element'') satisfying certain compatibility properties (see,
e.g., \cite{Kas} or  \cite{K1}). Namely, we define the
Gaussian to be the following element of the suitable completion of the
Grothendieck ring $K(\Cal C)$:
$$
\gamma_{\Cal C}=\sum_i C_i \dim X_i \<X_i\>,
$$
where $X_i$ are simple objects in $\Cal C$, $C_i$ is the value of the
Casimir element in $X_i$ (in \cite{K1}, these numbers are denoted by
$\theta_i$), $\dim X_i$ is the $q$-dimension of $X_i$, and $\<X_i\>$
is the class of $X_i$ in the Grothendieck ring.  In particular, if we
take the category of representations of the quantum group $\Uh$
corresponding to the Cartan subalgebra $\h\subset \g$ considered as a
commutative Lie algebra, then its irreducible representations are
parametrized by $\l\in P$, and they are all one-dimensional. One can
check that defining the universal $R$-matrix by $R|_{X_\l\o
X_\mu}=q^{(\l, \mu)}$, and the Casimir element  by
$C|_{X_\l}=q^{\l^2}$ endows this category with a structure of ribbon
category.  Thus, the Gaussian $\gamma_{\Uh}=\gamma$ for this category
is exactly given by the formula (3). On the other hand, if we consider
the category of representations of the quantum group $\Uq$, then the
Casimir element $C$ in this category is defined by 
$C=q^{2\rho}u^{-1}$, where $u$ is the Drinfeld's element (see details
in \cite{Kas, Chapter XIV.6}, where the element $\theta=C^{-1}$ is
discussed), and $C|_{L_\l}=q^{(\l ,\l+2\rho)}$. Thus, 
Gaussian for this category is given by $\gamma_{\Uq}=\sum_{\nu\in
P_+}q^{(\nu, \nu+2\rho)}(\dimq L_\nu) \chi_\nu$.  Therefore,
Proposition~1 can be rewritten as
$$
\gamma_{\Uq}=\biggl(\prod_{\a\in R_+}\frac{1}{1-q^{2(\a, \rho)}}
	\biggr)\gamma_{\Uh},
$$
which is closely connected  with the  Weyl formula for
a compact group $G$, which relates the measure on $G/\text{Ad }G=T/W$
coming from the Haar measure on $G$ with the Haar measure on $T$.

\head 2. Proof of Cherednik--Macdonald--Mehta identities
\endhead
In this section, we give a proof of the Cherednik--Macdonald--Mehta
identities (1) for $\g=\sln$. The proof is based on the realization of
Macdonald's polynomials as ``vector-valued characters'' for the
quantum group $\Un$, which was given in \cite{EK1}. For the sake of
completeness, we briefly outline these results here, referring the
reader to the original paper for a detailed exposition.

Let us fix $k\in \Z_+$ and denote by $U$ the finite-dimensional
representation of $\Un$ with highest weight $n(k-1)\omega_1$, where
$\omega_1$ is the first fundamental weight. We identify the zero
weight subspace $U[0]$, which is one-dimensional, with $\Cq$. 

For $\l\in P_+$, we denote by $\Phi_\l$ the unique intertwiner 
$$
\Phi_\l:L_{\l+(k-1)\rho}\to L_{\l+(k-1)\rho }\o U
$$
and define the generalized character $\ph_\l\in \Cq[P]\o U[0]\simeq
 \Cq[P]$  by 
$\ph_\l(q^x)=\Tr_{L_{\l+(k-1)\rho}}(\Phi_\l q^x)$.

We can now summarize the results of \cite{EK1} as follows:
$$\gathered
\ph_0=\prod_{\a\in
R_+}\prod_{i=1}^{k-1}(e^{\a/2}-q^{-2i}e^{-\a/2})=\de_k/\de\\
\ph_\l/\ph_0=P_\l
\endaligned 
\tag 6
$$
where $P_\l$  is the Macdonald polynomial with parameters
$q^2, t=q^{2k}$. 

We can also rewrite Macdonald's inner product in terms of the
generalized characters as follows. Recall that  Macdonald's  inner product 
on $\Cq[P]$ is defined by 
$$
\<f, g\>_k=\frac{1}{|W|}\int \de_k\overline{\de_k}f\bar g \, dx
$$
(this differs by a certain power of $q$ from the original definition
of Macdonald).  Obviously, one has
$$
\<P_\l, P_\mu\>_k=\<\ph_\l, \ph_\mu\>_1.
$$
In order to rewrite this in terms of representation theory, let
$\omega$ be the Cartan involution in $\Un$ (see \cite{EK1}). For a
$\Un$-module $V$, we denote by $V^\omega$ the same vector space but
with the action of $\Un$ twisted by $\omega$. Note that for
finite-dimensional $V$, we have $V^\omega\simeq V^*$ (not
canonically). Similarly, for an intertwiner $\Phi:L\to L\o U$ we
denote by $\Phi^\omega$ the corresponding intertwiner $L^\omega\to
U^\omega\o L^\omega$. Finally, for $\Phi_1:L_1\to L_1\o U,
\Phi_2:L_2\to L_2\o U$, define $\Phi_1\odot \Phi_2^\omega\in
\End(L_1\o L_2^\omega)$
as the composition $L_1\o L_2^\omega\to L_1\o U\o U^\omega\o L_2^\omega\to
L_1\o L_2^\omega$, where the first arrow is given by $\Phi_\l\o
\Phi_\mu^\omega$, and the second by the invariant
pairing $U\o U^\omega\to \Cq$
(which is unique up to a constant).  Then it was shown in \cite{EK1} that 
$$
(\ph_\l\overline{\ph_\mu})(q^x)=\Tr_V((\Phi_\l\odot\Phi^\omega_\mu)
q^{\Delta(x)})
=\sum_{\nu \in P_+}\chi_\nu(q^x)C_{\l\mu}^\nu
$$
where $V=L_{\l+(k-1)\rho}\o L^\omega_{\mu+(k-1)\rho}$ and
$C_{\l\mu}^\nu$ is the trace of $\Phi_\l\odot\Phi_\mu$ acting in the
multiplicity space $\Hom(L_\nu, V)$. 
As a corollary, we get the following result:
$$
\frac{1}{|W|}\int \de \bar\de \ph_\l\overline{\ph_\mu}
\bigl(\sum_{\nu\in P^+}a_\nu\chi_\nu\bigr) \, dx
= \sum_{\nu\in P^+} a_{\nu^*} C_{\l\mu}^\nu,
\tag 7
$$
where $\nu^*=-w_0(\nu)$ is the highest weight of the module $(L_\nu)^*$
(here $w_0$ is the longest element of the Weyl
group).

Of course, the coefficients $C_{\l\mu}^\nu$ are very difficult to
calculate. However, the formula above is still useful. For example, it
immediately shows that $\<\ph_\l, \ph_\mu\>_1=0$ unless $\l=\mu$,
which was the major part of the proof of the formula
$\ph_\l/\ph_0=P_\l$ in \cite{EK1}. It turns out that this formula also
allows us to prove the Cherednik-Macdonald--Mehta identities.

\proclaim{Theorem 2} Let $\ph_\l$ be the renormalized Macdonald
polynomials for the root system $A_{n-1}$ given by \rom{(6)}, and let
$\gamma$ be the Gaussian \rom{(3)}. Then
$$\aligned
\frac{1}{|W|} \int \de\overline{\de} \ph_\l \overline{\ph_\mu} \gamma\, dx
=& q^{(\l+k\rho)^2}q^{(\mu+k\rho)^2}\ph_\mu(q^{-2(\l+k\rho)})\\
&\quad\times\biggl(\prod_{\a\in R_+}(1-q^{2(\a, \rho)})
	\biggr)
q^{-2\rho^2}\|P_\l\|^2 \dimq L_{\l+(k-1)\rho},
\endaligned
\tag 8
$$
where $\|P_\l\|^2=\<P_\l, P_\l\>_k$. 
\endproclaim

\demo{Proof} From (7) and (4), we get 
$$
\int \de\overline{\de} \ph_\l \overline{\ph_\mu} \gamma\, dx= 
\biggl(\prod_{\a\in R_+}(1-q^{2(\a, \rho)})
	\biggr)\sum_{\nu\in P^+} q^{(\nu, \nu+2\rho)}(\dimq L_\nu)
	C_{\l\mu}^\nu.
\tag 9
$$

On the other hand, let $C$ be the Casimir element for $\Uq$ discussed
above.   Consider the intertwiner $(\Phi_\l\odot\Phi_\mu^\omega)
\Delta(C):V\to V$, where, as before, $V=L_{\l+(k-1)\rho}\o
L^\omega_{\mu+(k-1)\rho}$. Then it follows from $C|_{L_\l}=q^{(\l,
\l+2\rho)}$ that
$$
\Tr_V((\Phi_\l\odot\Phi_\mu^\omega)
\Delta(C) \Delta (q^{2\rho}))=\sum_{\nu\in P_+}C_{\l\mu}^\nu q^{(\nu,
\nu+2\rho)}\dimq L_\nu  
$$
which is exactly the sum in the right hand side of (9). On the other
hand, using 
$\Delta(C)=(C\o C)(R^{21}R)$, we can write 
$$
\Tr_V((\Phi_\l\odot\Phi_\mu^\omega)
\Delta(C) \Delta(q^{2\rho}))=q^{-2\rho^2}q^{(\l+k\rho)^2}q^{(\mu+k\rho)^2}
\Tr_V((\Phi_\l\odot\Phi_\mu^\omega) (R^{21}R) \Delta(q^{2\rho}))
$$
This last trace can be calculated, which was done in \cite{EK2,
Corollary~4.2}, and the answer is given by 
$$
\Tr_V((\Phi_\l\odot\Phi_\mu^\omega) (R^{21}R) \Delta(q^{2\rho}))
=\ph_\mu(q^{-2(\l+k\rho)})\|P_\l\|^2\dimq L_{\l+(k-1)\rho}.
$$
Combining these results, we get the statement of the theorem.
\qed\enddemo

The norms $\|P_\l\|^2$ appearing in the right-hand side of (8) are
given by Macdonald's inner
product identities 
$$
\|P_\l\|^2=\prod_{\a\in R_+}\prod_{i=1}^{k-1}
\frac {1-q^{-2 (\a, \l+k\rho)-2i}}
      {1-q^{-2 (\a, \l+k\rho)+2i}},
$$
which were conjectured in  \cite{M1, M2} and proved for root system
$A_{n-1}$ 
by Macdonald himself \cite{M3}; see also \cite{EK2} for the proof based
on representation theory of $\Un$, and \cite{Ch1} or a review in
\cite{K2} for a proof for arbitrary root systems.  Using this formula
and rewriting the statement of Theorem~2 in terms of Macdonald
polynomials $P_\l$ rather than $\ph_\l$, we get the
Cherednik--Macdonald--Mehta identities (1).

\remark{Remarks}
1. Note that the left-hand side of (8)   is symmetric in $\l,
\mu$. Thus, the same is true for the right-hand side, which is exactly the
statement of Macdonald's symmetry identity (compare with the proof in
\cite{EK2}).

2. The proof of Cherednik--Macdonald--Mehta identities given above
easily generalizes to the case when $q$ is a root of unity (see
\cite{K1} for the discussion of the appropriate
representation-theoretic setup). In this case, we need to replace the
set $P_+$ of all integral dominant weights by an appropriate (finite)
Weyl alcove $C$ (see \cite{K1}), and the integral $\int \de\bar\de f\, dx$
should be replaced by 
$\text{const}\sum_{\l\in C}f(q^{2(\l+\rho)})\dimq L_\l$.
Using the following obvious property of the Gaussian:
$$
\gamma(q^{2(\l+\rho)})=q^{-(\l, \l+2\rho)}\gamma(q^{2\rho})
$$
(which in this case coincides with formula (1.7)  in
\cite{K1}), 
it is easy to see that in this case the
Cherednik--Macdonald--Mehta identities are equivalent to
$$
S^{-1}T^{-1}S=TST
$$
where the matrices $S, T$  are defined in \cite{K1, Theorem~5.4}. This
identity is a part of a more general result, namely, that these
matrices $S, T$ give a projective representation of the modular group
$SL_2(\Z)$ on the space of generalized characters (see \cite{K1,
Theorem 1.10} and references therein). 
\endremark

\enddocument
\end